\documentclass[11pt,oneside]{amsart}

      \usepackage{amssymb}
			\usepackage{geometry}
			\usepackage{amsmath,amscd}
			\usepackage{hyperref}

      \theoremstyle{plain}
      \newtheorem{theorem}{Theorem}[section]

			\newtheorem*{theo}{Theorem}

      \theoremstyle{plain}
      \newtheorem{remarque}[theorem]{Remark}

      \makeatletter
      \def\@setcopyright{}
      \def\serieslogo@{}
      \makeatother
   
\begin{document}

%



   \author{Adrien Boyer}
   \address{Technion, Haifa, Israel}
   \email{adrienboyer@technion.ac.il}


  
   \title{Semisimple Lie groups satisfy property RD, a short proof}



   \subjclass{Primary 43; Secondary 22}


  
   \dedicatory{}



   \maketitle



\begin{abstract}
We give a short elementary proof of the fact that connected semisimple real Lie groups satisfy property RD. The proof is based on a process of linearization.

\end{abstract}

\section{Introduction}

A length function  $L:G \rightarrow \mathbb{R^{*}_{+}}$ on a locally compact group $G$ is a measurable function satisfying 

\begin{enumerate}
	\item $L(e)=0$ where $e$ is the neutral element of $G$
	\item $L(g^{-1})=L(g)$ 
	\item $L(gh)\leq L(g)+L(h).$
\end{enumerate}

A unitary representation $\pi:G \rightarrow U(H)$ on a complex Hilbert space has property RD with respect to $L$ if there exists $C>0$ and $d\geq 1$ such that for each pair of unit vectors $\xi$ and $\eta$ in $H$, we have $$\int_{G}\frac{|\left\langle \pi(g)\xi,\eta \right\rangle|^{2}}{(1+L(g))^{d}}dg\leq C$$ where $dg$ is a (left) Haar measure on $G$. We say that $G$ has property RD if its regular representation has property RD with respect to $L$. First established for free groups by U. Haagerup in \cite{H}, property RD
has been introduced and studied as such by P. Jolissaint in \cite{J}, who notably established it for groups of polynomial growth, and for classical hyperbolic groups. See \cite{V}[Chap. 8, p.69], and  for more details.\\

If $\pi$ denotes a unitary representation on a Hilbert space $H$, then $\overline{\pi}$ denotes its conjugate representation on the conjugate Hilbert space $\overline{H}$. The process of linearisation consists in working with $\sigma:G \rightarrow U(\overline{H} \otimes H)$ the unitary representation $\sigma=\overline{\pi} \otimes \pi$, see \cite{AB}[section 2.2].\\

A connected semisimple real Lie group with finite center can be written $G=KP$ where $K$ is a compact connected subgroup, and $P$ a closed amenable subgroup. We denote by $\Delta_{P}$ the right-modular function of $P$. Extend to $G$ the map $\Delta_{P}$ of $P$ as $\Delta : G \rightarrow \mathbb{R^{*}_{+}}$ with $\Delta(g)=\Delta(kp):=\Delta_{P}(p)$. It's well defined because $K\cap P$ is compact (observe that ${\Delta_{P}}_{|_{K\cap P}}=1$). The quotient $G/P$ carries a unique quasi-invariant measure $\mu$, such that the Radon-Nikodym derivative at $(g,x)\in G \times G/P$ denoted by  $c(g,x)=\frac{dg_{*}\mu}{d\mu}(x)$ with $g_{*}\mu(A)=\mu(g^{-1}A)$, satisfies $\frac{dg_{*}\mu}{d\mu}(x)=\frac{\Delta(gx)}{\Delta(x)}$ for all $g\in G$ and $x\in G/P$, (notice that for all $g\in G$, the function $x \in G/P \mapsto  \frac{\Delta(gx)}{\Delta(x)} \in \mathbb{R^{*}_{+}}$ is well defined). We refer to \cite{B}[Appendix B, Lemma B.1.3, p.~344-345] for more details. Consider the quasi-regular representation  $\lambda_{G/P}:G \rightarrow U(L^{2}(G/P))$ associated to $P$, defined by $(\lambda_{G/P}(g)\xi)(x)=c(g^{-1},x)^{\frac{1}{2}}\xi(g^{-1}x)$.  Denote by $dk$ the Haar measure on $K$, and under the identification $G/P=K/(K\cap P)$, denote by $d[k]$ the measure $\mu$ on $G/P$. \\

 The well-known Harish-Chandra function is defined by $\Xi(g):= \left\langle \lambda_{G/P}(g)1_{G/P},1_{G/P}\right\rangle$ where $1_{G/P}$ denotes the characteristic function of the space $G/P$.\\

In the rest of the paper we set $\sigma=\overline{\lambda_{G/P}}\otimes \lambda_{G/P}.$ Observe that $\overline{L^{2}(G/P)} \otimes L^{2}(G/P) \cong L^{2}(G/P \times G/P)$, via: $\xi \otimes \eta \longmapsto \left((x,y)\mapsto \overline{\xi(x)}\eta(y)\right)$. Notice that $\sigma$ preserves the cone of positive functions on $L^{2}(G/P \times G/P)$.\\

Let $G$ be a (non compact) connected semisimple real Lie group . Let $\frak{g}$ be its Lie algebra. Let $\theta$ be a Cartan involution. Define the bilinear form denoted by $(X,Y)$ such that for all $X,Y \in \frak{g}$, $(X,Y)=-B(X,\theta(Y))$ where $B$ is the Killing form. Set $|X|=\sqrt{(X,X)}$. 
Write $\frak{g}=\frak{l} \oplus \frak{p}$ the eigenvector space decomposition associated to $\theta$ ($\frak{l}$ for the eigenvalue $1$). Let $K$ be the compact subgroup defined as the connected subgroup whose Lie algebra $\frak{l}$ is the set of fixed points of $\theta$. Fix $\frak{a}\subset \frak{p}$ a maximal abelian subalgebra of $\frak{p}$. Consider the roots system $\Sigma$ associated to $\frak{a}$ and let $\Sigma^{+}$ be the set of positive roots, and define the corresponding positive Weyl chamber as
$$\frak{a}^{+}:=\left\{H\in \frak{a}, \alpha(H)>0, \forall \alpha \in \Sigma^{+}\right\}.$$

Let $A^{+}=Cl(\exp( \frak{a^{+}}))$, where $Cl$ denotes the closure of $\exp(\frak{a^{+}})$.  Consider the corresponding polar decomposition $KA^{+}K$. 
Then define the length function $$L(g)=L(k_{1}e^{H}k_{2}):=|H|$$  where $g=k_{1} e{H} k_{2}$ with $e^{H}\in A^{+}$. Notice that $L$ is $K$ bi-invariant. The disintegration of the Haar measure on $G$ according to the polar decomposition is $$dg=dkJ(H)dHdk$$ where $dk$ is the Haar measure on $K$, $dH$ the Lebesgue measure on $\frak{a}^{+}$, and $$J(H)=\prod_{\alpha\in \Sigma^{+}}\left(\frac{ e^{\alpha(H)}-e^{-\alpha(H)}}{2} \right)^{n_{\alpha}}$$ where $n_{\alpha}$ denotes the dimension of the root space associated to $\alpha$. See  \cite{K}[Chap.V, section 5, Proposition 5.28, p.141-142], \cite{GV}[Chap. 2, \S 2.2, p.65]  and \cite{GV}[Chap 2, Proposition 2.4.6, p.73] for more details. \\

The aim of this note is to give a short proof of the following known result (\cite{C},\cite{He}).
\\

\begin{theo} (C. Herz.)
Let $G$ be a connected real semisimple Lie group with finite center.
Then $G$ has property RD with respect to L.

\end{theo}

See \cite{C}[Proposition 5.5 and Lemma 6.3] for the case $G$ has infinite center.

\section{Proof}
\begin{proof}

We shall prove that the quasi-regular representation has property RD with respect to $L$ defined above. This implies that the regular representation has property RD with respect to $L$ by Lemma 2.3 in \cite{S}. 
Write $G=KP$ where $K$ is a compact subgroup and $P$ is a closed amenable subgroup of $G$.
It's sufficient to prove that there exists $d_{0} \geq 1$ and $C_{0}\geq 0$ such that $\int_{G}\frac{\left\langle \lambda_{G/P}(g)\xi,\xi\right\rangle^{2}}{(1+L(g))^{d_{0}}}dg<C_{0}$, for positive functions $\xi$, with $||\xi||=1$.\\

Take $\xi \in L^{2}(G/P)$ such that $\xi\geq 0$, and $||\xi||=1$. Define the function 

\begin{align*}
F:  G/P \times G/P &\rightarrow  \mathbb{R}_{+}\\
(x,y) & \mapsto \int_{K}\sigma(k)(\xi\otimes\xi) (x,y) dk. 
\end{align*}
 For all $(x,y) \in G/P \times G/P$, we have by the Cauchy-Schwarz inequality:
\begin{align*}
 \int_{K}\sigma(k)(\xi\otimes\xi) (x,y) dk &=\int_{K} \xi(k^{-1}x)\xi(k^{-1}y) dk\\
&\leq \left(\int_{K} \xi^{2}(k^{-1}x) dk \right)^{\frac{1}{2}}  \left( \int_{K}\xi^{2}(k^{-1}y) dk\right)^{\frac{1}{2}}. 
 \end{align*} 
Observe that the function $f:x \in G/P \mapsto \int_{K} \xi^{2}(k^{-1}x) dk \in \mathbb{R}_{+}$ is constant. Indeed, fix $x\in G/P$ and let $y$ in $G/P$. Write $y=hx$ for some $h\in K$ (as $K$ acts transitively on $G/P$ ). By invariance of the Haar measure we have $f(y)=\int_{K} \xi^{2}(k^{-1}y) dk=\int_{K} \xi^{2}(k^{-1}hx) dk= \int_{K} \xi^{2}(k^{-1}x) dk=f(x)$. 
If $e$ is the neutral element in $G$, we write $[e]\in G/P$. We have for all $x\in G/P$, $f(x)=f([e])$. \\
Hence, for all $x \in G/P$ we have
\begin{align*}
\int_{K} \xi^{2}(k^{-1}x) dk &=\int_{K} \xi^{2}(k^{-1}[e]) dk    \\
&=\int_{K/K \cap P} \xi^{2}([k^{-1}]) d[k]\\
&=||\xi||^{2}=1.
 \end{align*} 
Therefore $||F||_{\infty}:=\sup \left\{F(x,y), (x,y)\in G/P \times G/P \right\}\leq 1$. Hence $0\leq F \leq 1_{G/P \times G/P}$, where $1_{G/P \times G/P}$ denotes the characteristic function of  $G/P\times G/P$.\\

 Let $r$ be the number of indivible positive roots in $\frak{a}$. We know that there exists $C>0$ such that for all $H \in \frak{a}$ where $e^{H} \in A^{+}$ we have $$\Xi(e^{H})\leq Ce^{-\rho(H)}\left(1+L(e^{H})\right)^{r} $$ with $\rho=\frac{1}{2}\sum_{\alpha\in \Sigma^{+}}n_{\alpha} \alpha \in \frak{a}^{+}$, see  \cite{GV}[Chap 4, Theorem 4.6.4, p.161]. Hence for $d_{0}> \dim(\frak{a})+2r$, we have $$\int_{\frak{a}^{+}} \frac{ \Xi^{2}(e^{H})}{(1+L(e^{H}))^{d_{0}}} J(H)dH<\infty.$$

 We obtain for all $d \geq 0$ and for all positive functions $\xi$, with $||\xi||=1$:

\begin{align*}
\int_{G}\frac{\left\langle \lambda_{G/P}(g)\xi,\xi\right\rangle^{2}}{(1+L(g))^{d}}dg &=\int_{G} \frac{\left\langle \lambda_{G/P}(g)\xi,\xi\right\rangle\left\langle \lambda_{G/P}(g)\xi,\xi\right\rangle}{(1+L(g))^{d}}dg \\
&=\int_{G} \frac{\left\langle \sigma (g) \xi\otimes\xi,\xi \otimes \xi \right\rangle}{(1+L(g))^{d}}dg \\
                                                                                     &=\int_{K}\int_{\frak{a}^{+}} \int_{K} \frac{\left\langle \sigma (k_{1}e^{H}k_{2})\xi\otimes\xi,\xi \otimes \xi \right\rangle}{(1+L(k_{1}e^{H}k_{2}))^{d}} J(H)dk_{1}dHdk_{2}\\
																																										&=\int_{K}\int_{\frak{a}^{+}} \int_{K} \frac{\left\langle \sigma (e^{H})\sigma(k_{2})\left(\xi\otimes\xi\right),\sigma(k_{1}^{-1})\left(\xi \otimes \xi\right) \right\rangle}{(1+L(e^{H}))^{d}} J(H)dk_{1}dHdk_{2}\\
																																									&=\int_{\frak{a}^{+}} \frac{\left\langle \sigma (e^{H})\left(\int_{K}\sigma(k_{2})\left(\xi\otimes\xi\right) dk_{2}\right),\left(\int_{K}\sigma(k_{1}^{-1})\left(\xi \otimes \xi\right) dk_{1}\right) \right\rangle}{(1+L(e^{H}))^{d}} J(H)dH.\\
&=\int_{\frak{a}^{+}} \frac{\left\langle \sigma (e^{H})F,F \right\rangle}{(1+L(e^{H}))^{d}} J(H)dH \\
&\leq \int_{\frak{a}^{+}} \frac{\left\langle \sigma (e^{H}) 1_{G/P \times G/P},1_{G/P \times G/P} \right\rangle}{(1+L(e^{H}))^{d}} J(H)dH \\
&= \int_{\frak{a}^{+}} \frac{\left\langle \lambda_{G/P} (e^{H})1_{G/P},1_{G/P} \right\rangle^{2}}{(1+L(e^{H}))^{d}} J(H)dH .\\
&=\int_{\frak{a}^{+}} \frac{\Xi^{2}(e^{H})}{(1+L(e^{H}))^{d}} J(H)dH \\
\end{align*}
 
  Take $d_{0}>\dim(\frak{a})+2r$ and $C_{0}=\int_{\frak{a}^{+}} \frac{ \Xi^{2}(e^{H})}{(1+L(e^{H}))^{d_{0}}} J(H)dH$. We have fund $d_{0}\geq 1$ and $C_{0}>0$ such that for all positive functions $\xi$ in $L^{2}(G/P)$ with $||\xi||=1$, we have  $$ \int_{G}\frac{\left\langle \lambda_{G/P}(g)\xi,\xi\right\rangle^{2}}{(1+L(g))^{d_{0}}}dg \leq C_{0} $$  as required. 
\end{proof}
\begin{remarque}
The same approach applies to algebraic semisimple Lie groups over local fields. See \cite{AR}[section 1, (1.3)] and \cite{W}[Lemme II.1.5.]. 
\end{remarque}

\begin{remarque}
It's not hard to see that this approach shows that the representations of the principal series of $G$ (of class one, see \cite{GV}[ (3.1.12) p.103] ) satisfy also property RD.
\end{remarque}

\subsection*{Acknowledgements}
We would like to thank C. Pittet for helpful criticisms, and P. Delorme to point out the papers \cite{AR} and \cite{W} for the the case of semisimple Lie groups over local fields.

\end{document}